\documentclass[letterpaper, 10 pt, conference]{ieeeconf}
\IEEEoverridecommandlockouts                              % This command is only needed if 
\usepackage{mathptmx}
\usepackage[scaled=.92]{helvet}
\let\proof\relax
\let\endproof\relax
\usepackage{amsthm,amssymb}
\usepackage{mathrsfs}
\usepackage[colorlinks,citecolor=black,urlcolor=blue]{hyperref}  %%check
\usepackage{algorithm}
\usepackage{algorithm}
\usepackage{algorithmic}
\usepackage{xcolor}
\usepackage[tbtags]{amsmath}
\usepackage{graphicx}
\usepackage{epstopdf}
\usepackage{epsfig}
\usepackage[export]{adjustbox} %for large figures
\usepackage[caption=false, font=footnotesize]{subfig}
\usepackage{caption}
\newtheorem{theorem}{Theorem}[section]
\newtheorem{lemma}[theorem]{Lemma}

\newtheorem{corollary}[theorem]{Corollary} 
\theoremstyle{definition}

\theoremstyle{remark}

\newcommand*{\R}{\mathbb{R}}
\newcommand*{\Po}{\text{Prox}}
\newcommand*{\E}{\mathbb{E}}
\newcommand{\norm}[1]{\left\lVert#1\right\rVert}
\newcommand{\rnorm}[1]{\left.#1\right\rVert}
\newcommand{\lnorm}[1]{\left\lVert#1\right.}

\newcommand{\Iprod}[2]{\left\langle #1,#2\right\rangle}

\title{\LARGE \bf
A Proximal Zeroth-Order Algorithm for Nonconvex Nonsmooth Problems
}

\author{Ehsan Kazemi$^{1}$ and Liqiang Wang$^{1}$% <-this % stops a space
\thanks{$^{1}$Department of Computer Science, University of Central Florida, Orlando, FL 32816, USA
        }%
}

\allowdisplaybreaks
\begin{document}
\thispagestyle{empty}
\pagestyle{empty}

\maketitle
  %\and%
  %\author{\fnms{Peter}
   % \snm{Cholak}\ead[label=e2]{cholak.1@nd.edu}}%increase label by
  %	% one for each
  %	% author
  %\address{Department of Mathematics\\
   % University of Notre Dame\\
   % Notre Dame IN 46556\\
   % USA\\
   % \printead{e2} }
  % \and \author{\fnms{???} \snm{???}\ead[label=e3]{???}}
  % \address{\printead{e3}} \affiliation{???}

  %%%% INSERT EACH AUTHOR'S FIRST INITIAL AND SURNAME%%%%

\begin{abstract}
In this paper, we focus on solving an important class of nonconvex optimization problems
which includes many problems for example signal processing over a networked 
multi-agent system and distributed learning over networks. Motivated by many applications in which the local objective function is the sum of smooth  
but possibly nonconvex part, and non-smooth but convex part 
subject to a linear equality constraint, this paper proposes a proximal 
zeroth-order primal dual algorithm (PZO-PDA) that accounts for the information structure of the problem. 
This algorithm
only utilize the zeroth-order information (i.e., the functional values) of smooth functions, yet the flexibility 
is achieved for applications that only noisy information of the objective function is accessible, 
where classical methods cannot be applied.
We prove convergence and rate of convergence for PZO-PDA. Numerical experiments are provided to validate the theoretical results.
\end{abstract}

% main matter with bibliography goes here
\section{Introduction}
Consider the following optimization problem 
\begin{equation}\label{problem1}
\min_{x\in X} f(x)+h(x),\qquad \text{s.t.}\,\,Ax=b
\end{equation} 
where $x\in\R^N$, $A$ $\in {\mathbb{R}}^{M\times N}$, $b\in\R^M$, and $X\subseteq \R^N$ is a closed convex set; $f(x):\R^N\to\R$ is a continuous smooth but possibly nonconvex function;   $h(x):\R^N\to\R$ is a convex but possibly lower semi-continuous nonsmooth function.

Problem \eqref{problem1} is an interesting class that can be found in many application domains including statistical learning, compressed sensing and image processing.

Through applying the particular structures of the problem as described above, this paper desires to develop an efficient proximal zeroth-order gradient algorithm for \eqref{problem1}, which enjoys the
advantages of fast convergence rate and low computation cost.   

\subsection{Related Work}
Because of the large-scale nature, it is often  impractical to explore the second order information in the solution process. Therefore, concerning the provided information of the functions in question, the existing algorithms that solve \eqref{problem1} can be divided into two categories: zeroth-order methods and first-order methods. The first-order methods acquire gradient information of the objective at a given point, where we assume that gradient information of augmented Lagrangian (AL) using namely a first-order oracle is available. A first-order AL based algorithm for nonconvex nonsmooth optimization has developed in \cite{burachik2010inexact}. In \cite{lan2016iteration} the iteration complexity for the AL method is analyzed for smooth and convex objective functions. Recently, a proximal algorithm (PG-EXTRA) is proposed in \cite{shi2015proximal}, which uses constant stepsize and achieves $o(1/R)$ rate for nondifferentiable but convex optimization. A recent research attention for solving problem \eqref{problem1} in \cite{hajinezhadperturbed} has been devoted to the
so-called perturbed proximal primal-dual algorithm (PProx-PDA) in which a primal gradient descent step is performed followed by an approximate dual gradient ascent step. Our algorithm is closely related to  PProx-PDA, which achieves a sublinear convergence to only an $\epsilon$-stationary point. On the other hand, recently, the alternating direction method of multipliers (ADMM)  has been widely used for solving nonsmooth optimization problems \cite{boyd2011distributed,chang2014multi}. ADMM needs restrictive assumptions on the problem types in order to achieve convergence, and only very recently, the convergence of ADMM for nonconvex problems is investigated \cite{zhang2010convergence}.  

Despite the efficiency of first-order methods on solving \eqref{problem1}, these methods can not deal with problems that the gradient information is not available and we can only get an estimation of function $f(x)$. According to the available informational structure of the objective functions, zeroth-order methods only acquire the objective (or component) function value at any point via the so-called stochastic zeroth-order oracle. 
Roughly speaking, one way to approximate the gradient of a function using these information from function $f$ is by calculating the difference of
the function on two points which are near enough and then divide
the difference by the distance between the two points.
In \cite{nesterov2011random}, Nesterov and Spokoiny  proposed a general zeroth-order based method and proved a convergence rate of $O(N/R)$ for a zeroth-order stochastic gradient method applied to nonsmooth convex problems. 

Based on \cite{nesterov2011random}, Ghadimi and Lan \cite{ghadimi2013stochastic} developed a stochastic zeroth-order gradient method for both convex and nonconvex problems and proved a convergence rate of $O(N/R)$ for a zeroth-order stochastic gradient method on nonconvex smooth problems. Duchi et al. \cite{duchi2015optimal} proposed a stochastic zeroth-order mirror descent based method for convex and strongly convex functions and proved an $O(1/\sqrt{R})$ rate for a zeroth-order stochastic gradient method on convex objectives. Gao et al. \cite{gao2014information} proposed a stochastic gradient ADMM method that allows only noisy estimations of function values to be accessible and they proved an $O(1/R)$ rate, under some requirements for smoothing parameter and batch size.  Recently, Lian et al. \cite{lian2016comprehensive} proposed an asynchronous stochastic optimization algorithm with zeroth-order methods and proved a convergence rate of $O(1/\sqrt{R})$. In \cite{hajinezhad2017zeroth}, a distributed zeroth-order optimization algorithm was proposed for nonconvex minimization. However, the algorithm
and analysis presented in that work are based on  bounding the successive dual variables with preceding primal variables which are not applicable to the problems such as \eqref{problem1} with a nonsmooth regularizer $h(x)$ and a general convex constraint.

Observe that the composite forms appear in various applications; examples include: 1) in a geometric median problem, $f$ is null function and $h$ is an $l_2$ norm term \cite{eiselt2011foundations}; 2) in sparse subspace estimation problem \cite{gu2014sparse}, $f$ is a linear function, and $h$ is a nonconvex regularization term that enforces sparsity \cite{zhang2010nearly}; 3) in consensus problem, the generic model can be formulated as 
\[
\min_{x} f(x)+h(x) := \sum_{i=1}^{N} \left(f_i(x_i)+h_i(x_i)\right), 
\]
where for each $i$ function $f_i$ represents nonconvex activation functions of neural networks \cite{allen2016variance}, or data fidelity term \cite{mateos2010distributed} such as squared $l_2$ norm. Function $h_i$ is a  regularization term such as $l_1$ norm or smooth $l_2$ norm \cite{boyd2011distributed}, or the indicator function for a closed convex set \cite{chang2014distributed}. 
\subsection{Contributions and Paper Organization}
We design an algorithm, belonging to zeroth-order  methods, to account for
the informational structure of the objective functions, which achieves optimality condition with provable global subblinear rate. In Section \ref{AlgoDevelopment}, we introduce a proximal zeroth-order primal dual algorithm for which we shall use PZO-PDA as its acronym. The proposed algorithm would allow parallel updates that can be used for linearly constrained multi-block structured optimization model of problem \eqref{problem1} for efficient parallel computing \cite{xu2018accelerated}.

Convergence and rate of convergence for PZO-PDA is established in Section \ref{converAnalysis}. We show that PZO-PDA converges sublinearly to stationary points of problem \eqref{problem1}. 
Numerically, the performance of PZO-PDA is demonstrated in Section \ref{Experiments}.
The results confirm that PZO-PDA performs efficiently and stably. To our knowledge, our algorithm is {the first} proximal zeroth-order  primal dual algorithms for nonconvex nonsmooth constrained optimizations with convergence guarantee. 

\subsection{Notation}
We use $\norm{.}$ for Euclidean norm. For a given vector $v$, and matrix $W$, we define $\norm{v}_W^2 := v^TWv$. We let $\nabla_i f(x)$ to denote the partial gradient of $f$ with respect to $x_i$ at $x$ and $\nabla f(x)$ with respect to $x$. For matrix $W$, $W^T$ represent its transpose. For two
vectors $a, b$, we use $\Iprod{a}{b}$ to represent their inner product. We let ${I_N}$ denote the identity matrix of size $N$. The indicator function for convex set $X$ is indicated by $\iota_X$ which is defined as $\iota_X(y)=1$ when $y\in X$, and $\iota_X(y)=0$ otherwise. For a nonsmooth convex function $h(x), \partial h(x)$ denotes the subdifferential set defined by
\begin{equation}
\partial h(x) = \{v\in\R^N; h(x)\geq h(y)+\Iprod{v}{x-y}\,\, \forall y\in \R^N\}.
\end{equation}
For a convex function $h(x)$ and a constant $\alpha > 0$ the proximity operator is defined as below
\begin{equation}
\displaystyle\text{prox}_{h}^{\alpha}(x) := \text{argmin}_{z}\,\, \frac{\alpha}{2}\norm{x-z}^2+h(z).
\end{equation}
\section{Algorithm Development}\label{AlgoDevelopment}
We assume for any given $x \in \text{dom}(f)$, we get a noisy approximation of the true function value $f(x)$ by calling a stochastic zeroth-order oracle $\mathcal{SZO}$, which returns a quantity denoted by $\mathcal{F}(x,\xi)$ with $\xi \in \R$ being a random variable. One intuitive way to use $\mathcal{SZO}$ is when for example only the input and output of deep neural networks (DNNs) are observable. This is the case, for instance, when the general structure of network is not known and the gradient computation via back propagation is forbidden; however, we can train the model by observation which avoids the need for learning substitute models \cite{chen2017zoo}. PZO-PDA can only get a noisy estimation of function value $f(x)$ by calling $\mathcal{SZO}$ which returns $\mathcal{F}(x,\xi)$.

Now since we can access the $\mathcal{SZO}$, we shall present some basic concepts proposed in \cite{nesterov2011random}, to approximate the first- order information of a given function $f$. Let $S$ be a unit ball in $\R^N$ and $U$ be the uniform distribution on $S$. Given $\mu > 0$, then the smoothing function $f_{\mu}$ is defined as  

\begin{equation}
f_{\mu}(z) = \E_{\{v\sim U\}}[f(z+\mu v)] = \frac{1}{\alpha(N)}\int_{S}{f(x+\mu v) dv} 
\end{equation}
where $\alpha(N)$ is the volume of $S$. Some properties of the smoothing function are shown in \cite{gao2014information}. If $f\in {\mathcal{C}}_L^1$, then $f_{\mu}\in \mathcal{C}_{L_{\mu}}^1$ with $L_{\mu} \leq L$ and
\begin{equation}\label{smoothed-func}
\begin{split}
\nabla f_{\mu}(x) &= \E_{\{v\sim U\}}\left[\frac{N}{\mu}f(x+\mu v)v\right]\\
& = \frac{N}{\beta(N)}\int_{v\in S}{\frac{f(x+\mu v) - f(x)}{\mu}}vdv.
\end{split} 
\end{equation}
where $\beta(N)$ is the surface area of unit ball $S$ in $\R^N$. In addition, for any $x\in \text{dom}(f)$, we have
\begin{equation}\label{smoothed-func-prop}
\begin{split}
&|f_\mu(x)-f(x)|\leq \frac{{L}\mu^2}{2}\\
&\norm{\nabla f_\mu(x)-\nabla f(x)}\leq \frac{\mu N{L}}{2}.
\end{split}
\end{equation}
Specifically, inspired by equation \eqref{smoothed-func}, by calling $\mathcal{SZO}$ which returns a quantity $\mathcal{F}(x,\xi)$, we define the zeroth-order stochastic gradient of $f$ 
\begin{equation}\label{zeroth-order-def}
G_{\mu}(x,v,\xi) = \frac{N}{\mu}\left({\mathcal{F}(x+\mu v,\xi)-\mathcal{F}(x,\xi)}\right) v
\end{equation}
where the constant $\mu >0$ is smoothing parameter and $v\in S$ is a uniform random vector.

Towards finding a solution for \eqref{problem1}, by utilizing the above notations and definitions, we present a zeroth-order primal-dual based scheme. Let us introduce the augmented Lagrangian (AL) function for problem \eqref{problem1} as
\begin{equation}
L_{\rho}\left(x,\lambda\right) = f(x)+h(x)+\left<\lambda,Ax-b\right>+\frac{\rho}{2}\norm{Ax-b}^2
\end{equation}  
where $\lambda\in\R^M$ is the dual variable, and $\rho > 0$ is the penalty parameter. By adding the proximal term $\frac{\beta}{2}\norm{x-y}^2_{B^TB}$, with $B\in \R^{M\times N}$, we introduce the following zeroth-order approximation of AL,  
\begin{equation}
\begin{split}
\widetilde{L}_{\rho,\gamma}\left(x,x^r,\lambda\right)&= u(x,x^r)+h(x)+\left<(1-\rho\gamma)\lambda,Ax-b\right>\\
&~~~~+\frac{\rho}{2}\norm{Ax-b}^2+\frac{\beta}{2}\norm{x-x^r}^2_{B^T\,B}
\end{split}
\end{equation}
where we introduce a new parameter $\gamma >0$ such that $\rho\gamma <1$. In above notation, we define surrogate function ${u(x,x^r) :=\left\langle \overline{G}_{\mu}^r,x-x^r\right\rangle}$ to be the linear approximation of $f(x)$, where $\overline{G}_{\mu}^{r}$ is given by 
\[
\overline{G}_{\mu}^r =\frac{1}{J_r}\sum_{j=1}^{J_r}{{G}_{\mu}(x^r,v^r_j,\xi^r_j)}
\]
and we set $v^r:= \{v_{j}^r\}_{j=1}^{J_r}$, $\xi^r:= \{\xi_{j}^r\}_{j=1}^{J_r}$.

The steps of the proposed PZO-PDA algorithm is described below (Algorithm \ref{alg:PProx-PDA}). 

%\begin{algorithm}
%\caption{The perturbed proximal primal-dual algorithm (PProx-PDA)}\label{alg:PProx-PDA}
%\begin{algorithmic}[1]
%\Initialize{$\lambda^0$ and $x^0$}
%\Repeat\,\, {update variables by}{}
%\State
%{\color{red} Generate $\phi^r_{j}\in \R^N, j=1,2,,\ldots, J$ from an i.i.d standard Gaussian distribution and calculate $\bar{G}_{\mu}(x^r,\phi^r,\xi^r)\in\R^N$ by
%\begin{equation*}
%\bar{G}_{\mu}(x^r,\phi^r,\xi^r) = \frac{1}{J}\sum_{j=1}^J\frac{H(x^r+\mu\phi_{j}^r,\xi_{j}^r)-H_i(x^r,\xi_{j}^r)}{\mu}\phi_{j}^r,
%\end{equation*}
%where we have defined $\phi^r:= \{\phi_{j}^r\}_{j=1}^J$, $\xi^r:= \{\xi_{j}^r\}_{j=1}^J$; Set  
%\begin{equation*}
%\begin{split}
%&u(x,x^r) :=\left\langle\bar{G}_{\mu}(x^r,\phi^r,\xi^r),x-x^r\right\rangle
%\end{split}
%\end{equation*}
%
%}
%\State\begin{equation}\label{update:primal}
%\begin{split}
%x^{r+1}&=\displaystyle\text{argmin}_{x\in X} u(x,x^r)+h(x)+\left<(1-\rho\gamma)\lambda^r,Ax-b\right>\\
%&~~+\frac{\rho}{2}\norm{Ax-b}^2+\frac{\beta}{2}\norm{x-x^r}^2_{B^T\,B}
%\end{split}
%\end{equation}
%\State \begin{equation}\label{update:dual}\lambda^{r+1} = (1-\rho\lambda)\lambda^r+\rho \left(Ax^{r+1}-b\right)
%\end{equation}
%\Until{\bf convergence} 
%\State {{{\bf Output}}} Iterate $x_a$ chosen uniformly random from $\{x^r\}_{r=0}^{T}$.
%\end{algorithmic}
%\end{algorithm}
\begin{algorithm}[H]
\caption{The proximal zeroth-order primal-dual algorithm (PZO-PDA)}\label{alg:PProx-PDA}
\begin{algorithmic}[H1]
\renewcommand{\algorithmicrequire}{\textbf{Input:}}
 \renewcommand{\algorithmicensure}{\textbf{Output:}}
\REQUIRE{ $x^0$}, $\lambda^0$, $R$, $\rho$, $\gamma$, $\beta$, $\mu$, $\{J_r\}_{r=1,\ldots,R}$ 
\FOR{$r=1$ \textbf{to} $R$ }

\STATE Generate $v^r_{j}\in \R^N, j=1,2,,\ldots, J_r$ from an i.i.d uniform distribution from unit ball in $\R^N$.
 
 \STATE At the $r$th iteration, we call $\mathcal{SZO}$, $J_r$ times to obtain ${G}_{\mu}(x^r,v^r_{j},\xi^r_{j})$, $j=1,\ldots,J_r$ by
\[
{G}_{\mu}(x^r,v^r_j,\xi^r_j) = \frac{N}{\mu}{\left(\mathcal{F}(x^r+\mu v_{j}^r,\xi_{j}^r)-\mathcal{F}(x^r,\xi_{j}^r)\right)}v^r_j.
\]

\STATE Then set $\overline{G}_{\mu}^r =\frac{1}{J_r}\sum_{j=1}^{J_r}{{G}_{\mu}(x^r,v^r_j,\xi^r_j)}$.

 \STATE Set  $u(x,x^r) :=\Iprod{\overline{G}_{\mu}^r}{x-x^r}$
, and compute

\begin{equation}\label{update:primal}
\begin{split}
x^{r+1}&=\displaystyle\text{argmin}_{x\in X} u(x,x^r)+h(x)+\left<(1-\rho\gamma)\lambda^r,Ax-b\right>\\
&~+\frac{\rho}{2}\norm{Ax-b}^2+\frac{\beta}{2}\norm{x-x^r}^2_{B^T\,B};
\end{split}
\raisetag{15pt}
\end{equation}

$\lambda^{r+1} = (1-\rho\gamma)\lambda^r+\rho \left(Ax^{r+1}-b\right)$.

\ENDFOR
\ENSURE  Iterate $x_a$ chosen uniformly random from $\{x^r\}_{r=0}^{R}$.
 \end{algorithmic}
\end{algorithm}
Algorithm \ref{alg:PProx-PDA} is related to the proximal method of multipliers first developed in \cite{rockafellar1976augmented}, however the theoretical results derived are only developed for convex problems. Each iteration of PZO-PDA performs a gradient descent step on the approximation of AL function, followed by taking one step of approximate dual gradient ascent. In fact, the primal variable $x$ is updated by minimization of function $\widetilde{L}$ rather
than AL function $L$, and it is posed as a problem to determine the convergence of iteration of this modified algorithm to the stationary solutions. The use of the surrogate function $u(x, x^r) := \Iprod{\overline{G}_{\mu}^r}{x-x^r}$ ensures
that only zeroth-order information is used for the primal update, where $\overline{G}_{\mu}^r$ is calculated by calling $\mathcal{SZO}$ multiple times at each iteration. 

Some remarks are in order here. The PZO-PDA is tightly associated to the classical Uzawa primal-dual method \cite{uzawa1958iterative}, which has been exploited to solve convex saddle point problems and linearly constrained convex problems \cite{nedic2009subgradient}. However, the primal and dual parameters are perturbed in approximation of AL function to facilitate convergence analysis. The appropriate choice of scaling matrix $B$ ensures problem \eqref{update:primal} is strongly convex.
In fact, matrix $B$ is often
used to eliminate the nonconvexity in the augmented Lagrangian, in order for the obtained subproblem to be strongly convex, or even to provide a closed-form solution through choosing matrix $B$ with $A^TA + B^TB \succeq I_N$. Although parameters $\rho$, $\gamma$ and $\beta$ are fixed for all $r$, it could be shown that adapting the parameters can accelerate the convergence of the algorithm \cite{hajinezhadperturbed}. Finally, we note that step 2 in Algorithm \ref{alg:PProx-PDA} is decomposable over the variables, therefore they are well-situated to be implemented in a distributed manner.

Before conducting the convergence analysis for  Algorithm \ref{alg:PProx-PDA}, let us first make some assumptions on $f(x)$ and $F(x,\xi)$. Functions $f(x)$ and $F(x,\xi)$, which is a noisy estimation of $f$ at $x$ when $\mathcal{ZSO}$ is called, satisfy 
\begin{enumerate}
\item [A1.] $\E_{\xi}[\mathcal{F}(x,\xi)] = f(x)$ and $\E_{\xi}[\nabla \mathcal{F}(x,\xi)] = \nabla f(x)$.
\item [A2.] The constant $\sigma \geq 0$ satisfies 
\[\E[ \norm{\nabla \mathcal{F}(x; \xi)-\nabla f(x)}^2]\leq \sigma^2.\]
\item [A3.] There exists $K\geq 0$ such we have $\norm{\nabla f(x)}\leq K$.
\end{enumerate}

Next we present some properties of function $G_{\mu}(x^r,v^r,\xi^r)$ defined in \eqref{zeroth-order-def}.
\begin{lemma}\label{free-grad-approx}\cite{gao2014information}
Suppose that $G_{\mu}(x^r,v^r,\xi^r)$ is defined as in \eqref{zeroth-order-def}, and assumptions A.1 and A.2 hold. Then
\begin{equation}
\E_{\xi^r,v}[G_{\mu}(x^r,v^r,\xi^r)] = \nabla f_{\mu}(x^r).
\end{equation}
If further assumption A.3 holds, then we have the following 
\begin{equation}
\E_{\xi^r,v}[\norm{\overline{G}_{\mu}^r-\nabla f_{\mu}(x^r)}^2]\leq \frac{{\tilde{\sigma}^2}}{J_r}
\end{equation}
where $\tilde{\sigma}^2:=2N[K^2+\sigma^2+\mu^2{L}^2N]$.
\end{lemma}
\section{Convergence Analysis}\label{converAnalysis}
In this part we analyze the behavior of PZO-PDA algorithm. Our analysis combines ideas from classical proof in \cite{gao2014information}, as well as two recent constructions \cite{hajinezhadperturbed,hong2017prox}. Our construction differs from the previous works in a number of ways, in particular, the constructed algorithms involved first-order methods, but are not applicable when only zeroth-order information of the objective function is available. Moreover the analysis in \cite{hajinezhadperturbed} only guarantees global convergence of $\{(x^r,\lambda^r)\}_{r\in\mathbb{N}}$ to an $\epsilon$-stationary point, while we show PZO-PDA converges to a stationary solution of 
 \eqref{problem1}--provided the sequence of iterates is bounded. Further, we use the optimality gap to measure the quality of the solution, which makes the analysis more involved compared with the existing global error measures in \cite{ghadimi2013stochastic}. 

In the sequel, we will frequently use the following identity 
\begin{equation}\label{eq14}
\left<b,b-a\right> = \frac{1}{2}\left(\norm{b-a}^2+\norm{b}^2-\norm{a}^2\right).
\end{equation} 
Further, for simplicity we define $w^r := \left(x^{r+1}-x^r\right) - \left(x^{r}-x^{r-1}\right)$.
We may assume without loss of generality, $J_r = J$, for all $r$.
To establish convergence and rate of convergence for PZO-PDA, we assume, without loss of generality, that $f(x), h(x)\geq 0$ for all $x\in X$.  We  choose matrix $B$ in Algorithm \ref{alg:PProx-PDA} to satisfy $AA^T+BB^T\succeq I_M$  in order to ensure that the strong convexity of regularization term dominates the nonconvexity of function $f(x)$.

We first analyze the dynamics of dual variables with running one iteration of PZO-PDA.
\begin{lemma}\label{iterationdecay}
Under Assumptions A, for all $r\geq 0$, the iterates of PZO-PDA satisfy 
\begin{equation}\label{lem1:eq17}
\begin{split}
& \frac{1-\rho\gamma}{2\rho} \norm{\lambda^{r+1}-\lambda^{r}}^2+\frac{\beta}{2} \norm{x^{r+1}-x^{r}}_{B^TB}^2\\
& \leq \frac{1-\rho\gamma}{2\rho} \norm{\lambda^{r}-\lambda^{r-1}}^2+\frac{\beta}{2} \norm{x^{r}-x^{r-1}}_{B^TB}^2\\
&~~~~+{3}\frac{\tilde{\sigma}^2}{J}+\frac{3L_{\mu}^2}{2}\norm{x^{r}-x^{r-1}}^2\\
&~~~~+ \frac{1}{2}\norm{x^{r+1}-x^{r}}^2 -\gamma\norm{\lambda^{r+1}-\lambda^{r}}^2\qquad \forall r \geq 1. 
\end{split}
\end{equation}
\end{lemma}
\proof The optimality condition of $x^{r+1}$ in \eqref{update:primal} is given by 
\begin{equation}\label{lem1:eq18}
\begin{split}
&\left<\overline{G}_{\mu}^{r}  +A^T \lambda^r(1-\rho\gamma) + \rho A^T(Ax^{r+1}-b)\right.\\
\,\,&\left.~~+\beta B^TB(x^{r+1}-x^{r})+\eta^{r+1},x^{r+1}-x\right> \leq 0,\,\,\,\forall\,\,x\in X
\end{split}
\end{equation}
where $\eta^{r+1}\in\partial h(x^{r+1})$. By considering the optimality condition of the same equation \eqref{update:primal} for $x=x^r$, we get  
\begin{equation}\label{lem1:eq19}
\begin{split}
&\left<\overline{G}_{\mu}^{r-1} +A^T \lambda^{r-1}(1-\rho\gamma) + \rho A^T(Ax^{r}-b)\right.\\
\,\,&\left.~~~+\beta B^TB(x^{r}-x^{r-1})+\eta^{r},x^{r}-x\right>\leq 0,\,\,\,\forall\,\,x\in X
\end{split}
\end{equation}
where $\eta^{r}\in\partial h (x^r)$. We set $x=x^r$ and $x=x^{r+1}$ in equations \eqref{lem1:eq18} and \eqref{lem1:eq19}, respectively and adding the resulting inequalities. Applying the dual update in PZO-PDA yields
\begin{equation}\label{lem1:eq20}
\begin{split}
&\left<\overline{G}_{\mu}^{r} - \overline{G}_{\mu}^{r-1},x^{r+1}-x^{r}\right > +\left< A^T (\lambda^{r+1}-\lambda^{r}),x^{r+1}-x^{r}\right>\\
\,\,&~~~+\beta \left<B^TBw^r,x^{r+1}-x^{r}\right>\\
&~~~\leq \left<\eta^{r}-\eta^{r+1},x^{r+1}-x^{r}\right>\leq 0,\,\,\,\forall\,\,x\in X
\end{split}
\end{equation}
where the last inequality follows from convexity of $h$. Next we find upper bounds for the terms on the left hand of \eqref{lem1:eq20}. First, by using Cauchy-Schwarz inequality, and Lemma \ref{free-grad-approx} we have 
\begin{equation}\label{lem1:eq21}
\begin{split}
&\left<\overline{G}_{\mu}^{r-1}-\overline{G}_{\mu}^{r},x^{r+1}-x^{r}\right >\\
 &= \left<\overline{G}_{\mu}^{r-1}-\overline{G}_{\mu}^{r}+\nabla f_{\mu}(x^r) - \nabla f_{\mu}(x^{r})\right.\\
&~~~+\left.\nabla f_{\mu}(x^{r-1}) - \nabla f_{\mu}(x^{r-1}),x^{r+1}-x^{r}\right >\\
&  \leq \lnorm{\overline{G}_{\mu}^{r-1}-\overline{G}_{\mu}^{r}+\nabla f_{\mu}(x^r) - \nabla f_{\mu}(x^{r})}\\
&~~~+\rnorm{\nabla f_{\mu}(x^{r-1}) - \nabla f_{\mu}(x^{r-1})}\norm{x^{r+1}-x^{r}} \\
& \leq \frac{1}{2}\lnorm{\overline{G}_{\mu}^{r-1}-\overline{G}_{\mu}^{r}+\nabla f_{\mu}(x^r) - \nabla f_{\mu}(x^{r})}\\
&~~~+ \rnorm{\nabla f_{\mu}(x^{r-1}) - \nabla f_{\mu}(x^{r-1})}^2 + \frac{1}{2}\norm{x^{r+1}-x^{r}}^2 \\
&  \leq \frac{3}{2}\norm{\overline{G}_{\mu}^{r-1}-\nabla f_{\mu}(x^{r-1})}^2 + \frac{3}{2}\norm{\overline{G}_{\mu}^{r}-\nabla f_{\mu}(x^r)}^2\\
&~~~+\frac{3}{2} \norm{ \nabla f_{\mu}(x^{r}) - \nabla f_{\mu}(x^{r-1})}^2\\
&~~~+ \frac{1}{2}\norm{x^{r+1}-x^{r}}^2 \\
& \leq \frac{3}{2}\frac{\tilde{\sigma}^2}{J}+\frac{3}{2}\frac{\tilde{\sigma}^2}{J}+\frac{3L_{\mu}^2}{2}\norm{x^{r}-x^{r-1}}^2+ \frac{1}{2}\norm{x^{r+1}-x^{r}}^2 \\
&={3}\frac{\tilde{\sigma}^2}{J}+\frac{3L_{\mu}^2}{2}\norm{x^{r}-x^{r-1}}^2+ \frac{1}{2}\norm{x^{r+1}-x^{r}}^2.
\end{split}
\end{equation}
Further, proceeding as in \cite[Lemma 1]{hajinezhadperturbed} we have  
\begin{equation}\label{lem1:eq22}
\begin{split}
&\Iprod{A^T(\lambda^{r+1}-\lambda^{r})}{x^{r+1}-x^r} =\frac{1}{2}\left(\frac{1}{\rho}-\gamma\right)\left(\norm{\lambda^{r+1}-\lambda^{r}}^2\right.\\
&~~~~-\norm{\lambda^{r}-\lambda^{r-1}}^2+\left.\norm{(\lambda^{r+1}-\lambda^{r})-(\lambda^{r}-\lambda^{r-1})}^2 \right)\\
&~~~~ + \gamma \norm{\lambda^{r+1}-\lambda^{r}}^2.
\end{split}
\raisetag{15pt}
\end{equation}
For the last term in \eqref{lem1:eq20}, according to \eqref{eq14}, we obtain
\begin{gather}\label{lem1:eq23}
\begin{split}
\beta &\Iprod{B^T\,B w^{r}}{x^{r+1}-x^{r}}\\
 &=\frac{\beta}{2}\left(\norm{x^{r+1}-x^{r}}^2_{B^TB}- \norm{x^{r}-x^{r-1}}^2_{B^TB}+\norm{w^{r}}^2_{B^TB} \right)\\
&\geq\frac{\beta}{2} \left(\norm{x^{r+1}-x^{r}}^2_{B^TB}- \norm{x^{r}-x^{r-1}}^2_{B^TB}\right).
\end{split}
\raisetag{15pt}
\end{gather}
Plugging inequalities \eqref{lem1:eq21}-\eqref{lem1:eq23} in \eqref{lem1:eq20}, we obtain the desired result.
\endproof
Next we analyze the dynamics of primal iterations. To begin with, we construct the function $C(x,\lambda)$ as follow
\begin{equation}\label{eq24}
\begin{split}
C(x,\lambda)&:=f_{\mu}(x)+h(x)+\left<(1-\rho\gamma)\lambda,Ax-b-\gamma\lambda\right>\\
&+\frac{\rho}{2}\norm{Ax-b}^2
\end{split}
\end{equation}
where $f_{\mu}(x)$ denotes the smoothed version of function $f(x)$ defined in \eqref{smoothed-func}.
First we present some properties of the function $C(x,\lambda)$ over iterations.
\begin{lemma}\label{Titerationdecay}
Suppose that ${\beta > 3L_\mu+1}$ and $\rho \geq \beta$. Then for all $r\geq 0$ the iterates of PZO-PDA satisfy
\begin{equation}\label{lem2:25}
\begin{split}
C(x^{r+1}&,\lambda^{r+1}) + \frac{(1-\gamma\rho)\gamma}{2}\norm{\lambda^{r+1}}^2\\
\leq& C(x^{r},\lambda^{r})+\frac{(1-\gamma\rho)\gamma}{2}\norm{\lambda^r}^2\\
&+\left(\frac{(1-\rho\gamma)(2-\rho\gamma)}{2\rho}\right)\norm{\lambda^{r+1}-\lambda^{r}}^2 \\
&{-\left(\frac{\beta-3L_{\mu}-1}{2}\right)\norm{x^{r+1}-x^{r}}^2+\frac{1}{2}\frac{\tilde{\sigma}^2}{J}.}
\end{split} 
\end{equation}
\end{lemma}
Let us construct the following potential function $Q_c$, parametrized by a constant $c>0$
\begin{equation}
\begin{split}
&Q_c(x^{r+1},\lambda^{r+1};x^{r},\lambda^{r}) 
:= C(x^{r+1},\lambda^{r+1}) + \frac{(1-\rho\gamma)\gamma}{2}\norm{\lambda^{r+1}}^2\\
& ~~~~+\frac{c}{2}\left(\frac{(1-\rho\gamma)}{\rho}\norm{\lambda^{r+1}-\lambda^{r}}^2\right.\\
&\left.~~~~+\beta \norm{x^{r+1}-x^{r}}^2_{B^TB}+3L_{\mu}^2\norm{x^{r+1}-x^{r}}^2\right).
\raisetag{15pt}
\end{split}
\end{equation}
%Then according to the previous two lemmas, one can conclude there are constants $a_1$, $a_2$, such that 
%\begin{equation}\label{potential:lowbd}
%\begin{split}
%&{\color{blue} Q_c(x^{r+1},\lambda^{r+1};x^{r},\lambda^{r}) - Q_c(x^{r},\lambda^{r};x^{r-1},\lambda^{r-1})}\\
%&{\color{blue} \leq -a_1\norm{\lambda^{r+1}-\lambda^{r}}^2 -a_2\norm{x^{r+1}-x^{r}}^2,}
%\end{split}
%\end{equation} 
%where ${\color{blue} a_1=\left((1-\rho\gamma)\frac{\gamma}{2}+c\gamma-\frac{1-\rho\gamma}{\rho}\right)}$, and ${\color{blue} a_2=\left(\frac{\beta-3L}{2}-c\,L\right)}$.
%Therefore, it is easy to observe that in order to make the $Q_c$ function descent, it is sufficient to have 
%\begin{equation}\label{alg:param}
%{\color{blue}(1-\rho\gamma)\frac{\gamma}{2} + c\gamma - \frac{1-\rho\gamma}{\rho} >0,\,\,\text{and}\, \beta > (3+2c)L.}
%\end{equation} 
We skip the subscript $c$ if $c = 1$. In the following we show that when the algorithm parameters are chosen properly, the potential function will decrease along the iterations. 
\begin{lemma}\label{lemma:potential:lowbd}
Suppose the assumptions made in Lemma \ref{Titerationdecay} are satisfied and additionally the parameters $\beta$, $\rho$ and $\gamma$ satisfy the following conditions,
\begin{equation}\label{alg:param}
\begin{split}
&\frac{(1-\rho\gamma){\gamma}}{2} + \gamma - \frac{1-\rho\gamma}{\rho} >0\\
~~~~~&\beta > (3+3 L_\mu)L_\mu+2.
\end{split}
\end{equation} 
Then we have the following
\begin{equation}\label{potential:lowbd}
\begin{split}
&Q(x^{r+1},\lambda^{r+1};x^{r},\lambda^{r}) - Q(x^{r},\lambda^{r};x^{r-1},\lambda^{r-1})\\
& ~~~= -a_1\norm{\lambda^{r+1}-\lambda^{r}}^2 -a_2\norm{x^{r+1}-x^{r}}^2 + \frac{7}{2}\frac{\tilde{\sigma}^2}{J}
\end{split}
\end{equation}
with \small{$a_1=\big(\frac{(1-\rho\gamma){\gamma}}{2}+\gamma-\frac{1-\rho\gamma}{\rho}\big)$} and \small{$a_2=\big(\frac{\beta-1}{2}-\frac{3 L_{\mu}^2}{2}-\frac{3L_{\mu}}{2}-\frac{1}{2}\big)$}.
\end{lemma}
The proofs of Lemmas \ref{Titerationdecay} and \ref{lemma:potential:lowbd} are postponed to Appendix. Next, we show the lower boundedness of potential function. To precisely state the convergence of Algorithm \ref{alg:PProx-PDA}, let us define 
\[
\widetilde{Q}(x^{r+1},\lambda^{r+1};x^{r},\lambda^{r}) = {Q}(x^{r+1},\lambda^{r+1};x^{r},\lambda^{r})- (r+1)\frac{7}{2} \frac{\tilde{\sigma}^2}{J}.
\]
By Lemma \ref{lemma:potential:lowbd}, function $\widetilde{Q}$ decreases at each iteration of PZO-PDA.
\begin{lemma}\label{lowerbdd-potentialfunc}
Suppose Assumptions A are satisfied, and the algorithm parameters are chosen according to \eqref{alg:param}. Then, 
\begin{equation}\label{bdforlambda}
\frac{\gamma(1-\rho\gamma)}{2}\norm{\lambda^{r+1}}^2 \leq {Q}^{0}+\frac{7r}{2} \frac{\tilde{\sigma}^2}{J}.
\end{equation}
Furthermore, given a fixed iteration number $R$, if the number of calls to $\mathcal{SZO}$ at each iteration is $J = R^2$, for some constant $\underline{Q}$ the iterates generated by PZO-PDA satisfy
\begin{equation}
 Q(x^{r+1},\lambda^{r+1};x^{r},\lambda^{r}) \geq \underline{Q} > -\infty, \qquad \forall r\geq 0.
\end{equation}  
\end{lemma}
\proof
By an induction argument and using the fact that the function $\widetilde{Q}$ is nonincreasing, \eqref{bdforlambda} can be proved for all $r$. 

Second, following similar analysis steps
presented in \cite{hajinezhadperturbed}, taking a sum over $R$ iterations of $C(x^{r+1},\lambda^{r+1})$ we obtain 
\begin{equation}\label{eq34}
\begin{split}
\sum_{r=1}^R C(x^{r+1},&\lambda^{r+1})\\
 &\geq \sum_{r=1}^R \left(f_{\mu}(x^{r+1})+h(x^{r+1})+\frac{\rho}{2}\norm{Ax^{r+1}-b}^2\right)\\
&~~~~~~+\frac{(1-\rho\gamma)^2}{2\rho}(\norm{\lambda^{R+1}}^2-\norm{\lambda^{1}}^2)\\
&\geq -\frac{(1-\rho\gamma)^2}{2\rho}\norm{\lambda^{1}}^2
\raisetag{15pt}
\end{split}
\end{equation}
where the last inequality comes from the fact that both $f_{\mu}$ and $h$ are lower bounded by $0$. Therefore, the sum of the $C(\cdot,
\cdot)$ function is lower bounded. From \eqref{eq34} and by selecting $J=R^2$, we conclude that $\sum_{r=1}^R \widetilde{Q}(x^{r+1},\lambda^{r+1};x^{r},\lambda^{r})$ is also lower bounded by $-\frac{(1-\rho\gamma)^2}{2\rho}\norm{\lambda^{1}}^2-7{\tilde{\sigma}^2}$ for any $R$. Since $\widetilde{Q}$ is nonincreasing, we conclude that the potential function $\widetilde{Q}$ is lower bounded by some constant $\underline{Q}$. 
%that is we have
%\begin{equation}
%\underline{P} \geq -\frac{(1-\rho\gamma)^2}{2\rho}\norm{\lambda^1}^2.
%\end{equation}
Thus, by the definition of function $\widetilde{Q}$, the potential function $Q$ is also lower bounded.
This completes the proof.
\endproof

In the rest of section, we let $\omega^r$ to denote $\{(v^i,\xi^i)\}_{i=1}^r$. Next we define the optimality gap that measures the progress of the algorithm and solution quality 
\begin{equation}
\begin{split}
{\Psi}^r := \Bigg\lVert x^r-&\text{prox}_{h+\iota_X}^{\beta} \left[x^{r}-\frac{1}{\beta}[f{(x^r)}+A^T \lambda^{r}]\right]\Bigg\lVert^2\\
&+\frac{1}{\beta^2}\norm{x^{r+1}-x^{r}}^2+\frac{1}{\rho^2}\norm{\lambda^{r+1}-\lambda^{r}}^2.
\end{split}
\end{equation}
Now we present the main convergence result which provides a rate of convergence for PZO-PDA. 
\begin{theorem}\label{theo-conv}
Suppose Assumptions A hold, and $a$ is uniformly sampled from $\{1,2,\ldots,R\}$. We assume the parameters are chosen according to \eqref{alg:param}. Then we have the following bound for optimality gap in expectation
\begin{equation}
\E[\Psi^a] \leq \frac{2V}{R}\E[Q^0-Q^R]+ 2\tilde{b}\frac{\hat{\sigma}^2}{J} + \frac{\mu^2{L}^2N^2}{2}.
\end{equation}
Moreover, by choosing $\gamma={\mathcal{O}}(\frac{1}{R})$ such that $\alpha_0:=\rho\gamma < 1$ is remained constant, we have \begin{equation*}
\begin{split}
\E[\Psi^a] +\norm{Ax^{a+1}-b}^2&\leq \frac{4V}{R}\E[Q^0-Q^R]\\
&~~~+ \frac{2Q^0}{R(1-\alpha_0)}+\tilde{c}\frac{\hat{\sigma}^2}{J} + {\mu^2{L}^2N^2}.
\end{split}
\end{equation*}
Here, $V$, $\tilde{b}$ and $\tilde{c}$ are constants which do not depend on the problem accuracy.
\end{theorem}
\proof
First, we let ${\Psi}^r_{\mu}$ denote the smoothed version of optimality gap which is defined as follows 
\begin{equation}
\begin{split}
{\Psi}^r_{\mu} := \Bigg\lVert x^r-&\text{prox}_{h+\iota_X}^{\beta} \left[x^{r}-\frac{1}{\beta}[f_{\mu}{(x^r)}+A^T \lambda^{r}]\right]\Bigg\lVert^2\\
&+\frac{1}{\beta^2}\norm{x^{r+1}-x^{r}}^2+\frac{1}{\rho^2}\norm{\lambda^{r+1}-\lambda^{r}}^2.
\end{split}
\end{equation}
The $x$-subproblem in Algorithm \ref{alg:PProx-PDA} equivalently can be formulated as
\begin{equation*}
x^{r+1} = \text{prox}_{h+\iota_X}^{\beta} \Big[x^{r+1}-\frac{1}{\beta}[\overline{G}_{\mu}^{r}+A^T \lambda^{r+1}+\beta B^TB(x^{r+1}-x^{r})]\Big].
\end{equation*}
Using the above equality and definition of optimality gap, we have
\begin{equation*}
\begin{split}
&\E_{\omega^r}[{\Psi}^r_{\mu}] = \E_{\omega^r}\norm{x^r-\text{prox}_{h+\iota_X}^{\beta} \left[x^{r}-\frac{1}{\beta}[f_{\mu}{(x^r)}+A^T \lambda^{r}]\right]}^2\\
&~+\frac{1}{\beta^2}\E_{\omega^r}\norm{x^{r+1}-x^{r}}^2+\frac{1}{\rho^2}\E_{\omega^r}\norm{\lambda^{r+1}-\lambda^{r}}^2\\
&= \E_{\omega^r}{\lnorm{ x^r-\text{prox}_{h+\iota_X}^{\beta} \left[x^{r}-\frac{1}{\beta}[f_{\mu}{(x^r)}+A^T \lambda^{r}]\right]-x^{r+1}}}\\
&~+\text{prox}_{h+\iota_X}^{\beta} \Big[x^{r+1}-\frac{1}{\beta}[\overline{G}_{\mu}^{r}+A^T \lambda^{r+1}\rnorm{+\beta B^TB(x^{r+1}-x^{r})]\Big]}^2\\
&~+\frac{1}{\beta^2}\E_{\omega^r}\norm{x^{r+1}-x^{r}}^2+\frac{1}{\rho^2}\E_{\omega^r}\norm{\lambda^{r+1}-\lambda^{r}}^2\\
&\stackrel{a}{\leq} 2\E_{\omega^r}\norm{x^{r+1}-x^{r}}^2+\frac{6}{\beta^2}\E_{\omega^r}\norm{f_{\mu}{(x_r)}-\overline{G}_{\mu}^{r}}^2\\
&~~~~+\frac{6}{\beta^2}\E_{\omega^r}\norm{A^T \lambda^{r+1}-A^T \lambda^{r}}^2\\
&~~~~+6\E_{\omega^r}\norm{(I_N-B^TB)(x^{r+1}-x^{r})}^2\\
&~~~~+\frac{1}{\beta^2}\E_{\omega^r}\norm{x^{r+1}-x^{r}}^2+\frac{1}{\rho^2}\E_{\omega^r}\norm{\lambda^{r+1}-\lambda^{r}}^2\\
&\stackrel{b}{\leq} (2+6\sigma_{max}^2(\hat{B}^T\hat{B})+\frac{1}{\beta^2})\E_{\omega^r}\norm{x^{r+1}-x^{r}}^2\\
&~~~~+(\frac{6\sigma_{max}^2(A^TA)}{\beta^2}+\frac{1}{\rho^2})\E_{\omega^r}\norm{\lambda^{r+1}-\lambda^{r}}^2+ \frac{6\hat{\sigma}^2}{\beta^2 J}
\end{split}
\end{equation*}
where $\sigma_{max}$ denotes the largest eigenvalue of a matrix and we define $\hat{B} := I-B^TB$. In $\stackrel{a}{\leq}$ we applied the non-expansiveness of the proximal operator and in $\stackrel{b}{\leq}$ we used Lemma \ref{free-grad-approx}. Therefore,
\begin{equation}\label{eq56}
\begin{split}
\E_{\omega^r}[{\Psi}^r_{\mu}]\leq &b_1\E_{\omega^r}\norm{x^{r+1}-x^{r}}^2\\
&+b_2\E_{\omega^r}\norm{\lambda^{r+1}-\lambda^{r}}^2 + \frac{6\hat{\sigma}^2}{\beta^2 J}
\end{split}
\end{equation}
where $b_1=2+6\sigma_{max}^2(\hat{B}^T\hat{B})+\frac{1}{\beta^2}$ and $b_2=\frac{6\sigma_{max}^2(A^TA)}{\beta^2}+\frac{1}{\rho^2}$. Using \eqref{eq56} with the descent estimate for the potential function $Q$ in \eqref{potential:lowbd}, leads to
\begin{equation}\label{eq127}
\begin{split}
\E_{\omega^r}[{\Psi}^r_{\mu}] \leq V\E_{\omega^r}\big[&Q(x^{r},\lambda^{r};x^{r-1},\lambda^{r-1})\\
& - Q(x^{r+1},\lambda^{r+1};x^{r},\lambda^{r})\big]+ \tilde{b}\frac{\hat{\sigma}^2}{J}
\end{split}
\end{equation}
where we defined $V:= \frac{\max(b_1,b_2)}{\min(a_1,a_2)}$ and $\tilde{b} := -\frac{7}{2}V+\frac{6}{\beta^2}$. Summing \eqref{eq127} over $r=1$ to $R$, and divide both sides by $R$, we obtain
\begin{equation}\label{eq127-V1}
\begin{split}
\frac{1}{R}\sum_{r=1}^R&\E_{\omega^r}[{\Psi}^r_{\mu}]\\
\leq & \frac{V}{R}\E[Q(x^{1},\lambda^{1};x^{0},\lambda^{0}) - Q(x^{R+1},\lambda^{R+1};x^{R},\lambda^{R})]+ \tilde{b}\frac{\hat{\sigma}^2}{J}\\
=&  \frac{V}{R}\E[Q^0-Q^R]+ \tilde{b}\frac{\hat{\sigma^2}}{J}.
\end{split}
\raisetag{15pt}
\end{equation}
Now let us bound the gap $\Psi^r$. Using the definition of $\Psi^r$ we have 
\begin{equation*}
\begin{split}
&\E_{\omega^r}[\Psi^r] = \E_{\omega^r}\norm{x^r-\text{prox}_{h+\iota_X}^{\beta} \left[x^{r}-\frac{1}{\beta}[f{(x^r)}+A^T \lambda^{r}]\right]}^2\\
&~+\frac{1}{\beta^2}\E_{\omega^r}\norm{x^{r+1}-x^{r}}^2+\frac{1}{\rho^2}\E_{\omega^r}\norm{\lambda^{r+1}-\lambda^{r}}^2\\
&= \E_{\omega^r} {\left\lVert x^r-\text{prox}_{h+\iota_X}^{\beta} \left[x^{r}-\frac{1}{\beta}[f{(x^r)}+A^T \lambda^{r}]\right]\right.}\\
&~~~- \text{prox}_{h+\iota_X}^{\beta} \left[x^{r}-\frac{1}{\beta}[f_{\mu}{(x^r)}+A^T \lambda^{r}]\right]\\
&~~~+ \left.\text{prox}_{h+\iota_X}^{\beta} \left[x^{r}-\frac{1}{\beta}[f_{\mu}{(x^r)}+A^T \lambda^{r}]\right]\right\lVert\\
&~~~+\frac{1}{\beta^2}\E_{\omega^r}\norm{x^{r+1}-x^{r}}^2+\frac{1}{\rho^2}\E_{\omega^r}\norm{\lambda^{r+1}-\lambda^{r}}^2\\
&\leq 2\E_{\omega^r}[\Psi^r_{\mu}] + \frac{\mu^2{L}^2N^2}{2}
\end{split}
\end{equation*}
where the last inequality obtained from non-expansiveness of the proximal operator and \eqref{smoothed-func-prop}. Now sum over all iteration to obtain
\begin{equation}\label{eq128}
\begin{split}
\frac{1}{R}\sum_{r=1}^R\E_{\omega^r}[{\Psi}^r] \leq & 
\frac{2}{R}\sum_{r=1}^R\E_{\omega^r}[\Psi^r_{\mu}] + \frac{\mu^2{L}^2N^2}{2}\\
\leq&  \frac{2V}{R}\E[Q^0-Q^R]+ 2\tilde{b}\frac{\hat{\sigma}^2}{J} + \frac{\mu^2{L}^2N^2}{2}
\end{split}
\end{equation}
where in the last inequality we used \eqref{eq127-V1}.
Using the above inequality and the definition of
$x_a$ in Algorithm \ref{alg:PProx-PDA}, we obtain the desired result.

The second part follows by applying the dual update in PZO-PDA which gives
\begin{equation}
\begin{split}
\norm{Ax^{r+1}-b}^2 &\leq \frac{1}{\rho^2}\norm{\lambda^{r+1}-\lambda^{r}}^2 + \norm{\gamma\lambda^{r}}^2\\
&\leq {\Psi}^r+\frac{2Q^0\gamma}{(1-\rho\gamma)} + \frac{7(r-1)\gamma}{(1-\rho\gamma)} \frac{\tilde{\sigma}^2}{J}
\end{split}
\end{equation}
where the last inequality is obtained from definition of  ${\Psi}^r$  and Lemma \ref{lowerbdd-potentialfunc}. Summing up the above inequality for $r = 1, \ldots ,R$  and substituting \eqref{eq128}, by using $\gamma={\mathcal{O}}(\frac{1}{R})$ and $\alpha_0 = \rho\gamma$ we have
\begin{equation}\label{eq129}
\begin{split}
\frac{1}{R} \sum_{r=1}^R\norm{Ax^{r+1}-b}^2 &\leq \frac{2V}{R}\E[Q^0-Q^R]+ 2\tilde{b}\frac{\hat{\sigma}^2}{J} + \frac{\mu^2{L}^2N^2}{2}\\&~~~~+\frac{2Q^0}{R(1-\alpha_0)} + \frac{7}{(1-\alpha_0)} \frac{\tilde{\sigma}^2}{J}.\\
\end{split}
\raisetag{15pt}
\end{equation}
Finally, combining \eqref{eq128} and \eqref{eq129} yields the desired result.
\endproof
Note that while the first statement of the last theorem provides a rate of convergence for the optimality gap, the second part shows the size of the constraint violation also converges zero with the same order. 
In the following corollary, we comment on the structure of the proposed algorithm.
\begin{corollary}
Under the assumptions we made in Theorem \ref{theo-conv}, given a fixed iteration number $R$, if the smoothing parameter is chosen to be $\mu\leq{\frac{1}{\sqrt{R}}}$, and the number of calls to $\mathcal{SZO}$ at each iteration is $J = R$, then we have
\begin{equation}
\begin{split}
\E[\Psi^a] +\norm{Ax^{a+1}-b}^2\leq \frac{4V}{R}&\E[Q^0-Q^R]+ \frac{2Q^0}{R(1-\alpha_0)}\\
&+\tilde{c}\frac{\hat{\sigma}^2}{R} + \frac{{L}^2N^2}{R}.
\end{split}
\end{equation}
\end{corollary}
 Suppose that $\{(x^r,\lambda^r)\}_{r\in\mathbb{N}}$  is bounded and we let $(x^*,\lambda^*)$ denote any limit point of the sequence $\{(x^r,\lambda^r)\}_{r\in\mathbb{N}}$, then for a converging subsequence $(x^{r_j},\lambda^{r_j})\to (x^*,\lambda^*)$ we have  
\begin{equation}\label{remark-cons}
Ax^*-b = 0.
\end{equation}
Note that according to \eqref{lem1:eq18} the optimality condition of $x^{r+1}$ is given by 
\begin{equation*}
\begin{split}
&\left<\overline{G}_{\mu}^{r}  +A^T \lambda^r(1-\rho\gamma) + \rho A^T(Ax^{r+1}-b)\right.\\
\,\,~~&\left.~~~+\eta^{r+1}+\beta B^TB(x^{r+1}-x^{r}),x^{r+1}-x\right> \leq 0,\,\,\,\forall\,\,x\in X
\end{split}
\end{equation*}
where $\eta^{r+1}\in \partial h(x^{r+1})$. Therefore, from the above inequality combined with dual variable update overall we have 
\begin{equation}
\begin{split}
\left<\overline{G}_{\mu}^{r}\right.&+A^T \lambda^{r+1}+\eta^{r+1}\\
&\left. +\beta B^TB(x^{r+1}-x^{r}),x^{r+1}-x\right>\leq 0,\,\,\,\forall\,\,x\in X.
\end{split}
\end{equation}
If $\E[\overline{G}_{\mu}^{r}] \to \nabla f(x^*)$, this inequality and the limit $x^{r+1}-x^r\to 0$ imply the following optimality condition
\begin{equation}
 \Iprod{ \nabla f(x^*)+A^T \lambda^* + \eta^*}{x^*-x}\leq 0,\,\,\,\forall\,\, x\in X
\end{equation}
where $\eta^*$ is some vector that satisfies $\eta^* \in \partial h(x^*)$. The inequality above with \eqref{remark-cons} show $(x^*,\lambda^*)$ is a stationary point of problem \eqref{problem1}.

\section{Experimental Results}\label{Experiments}
Numerical experiments are performed over a connected network
consisting of {$N = 10$} agents {and 27 bidirectional edges}. In particular, we focus on
the constrained non-negative principle component analysis (PCA) problem, which can be formulated as $l_1$ regularized
least squares problem in the form
\[
\min_{x} \sum_{i=1}^n f_i(x_i) + \sum_{i=1}^n h_i(x_i),\qquad s.t \,Ax=0,\,\,\norm{x_i}^2\leq 1
\]
with $f_i(x_i)=-\Iprod{x_i}{Z_ix_i}$, $h_i(x)={\alpha}\norm{x_i}_1$, $x=\{x_i\}_{i=1}^N$ and $X$ in problem \eqref{problem1} is $X=\left\{x=\{x_i\}_{i=1}^N\left\vert\right.\norm{x_i}^2\leq 1, i=1,\ldots,N\right\}$.
Here, $\alpha$ is the regularization parameter on agent $i$ and $Z_i = M_i^T\,M_i\in\R^{d\times d}$, where $M_i\in\R^{p\times d}$ is the measurement matrix, and $p$ is the batch size. This problem have applications in decentralized multi-agent compressive sensing problem \cite{mardani2013decentralized}, where the goal of the agents is to jointly estimate the sparse signal $x$. 
In experiments, the sparse signal $x_i$ has dimension $d = 10$, yet each agent holds $p=100$ measurements, and its regularization parameter is {$\alpha = 10^{-4}$}. 

The elements of the measurement matrices $M_i$
are generated randomly following
uniform distribution in the interval $(0,1)$.
We set $\gamma = 10^{-5}$ and the penalty parameter $\rho$ and $\beta$ in Algorithm \ref{alg:PProx-PDA} are chosen to  fulfill theoretical bounds given in \eqref{alg:param}. The smoothing parameter is set $\mu = \frac{1}{\sqrt{R}}$, and the maximum number of iterations is chosen {$R = 500$}. The noise $\xi$ is generated from i.i.d Gaussian distribution with mean $0$ and standard deviation $0.01$. The initial solution of signal $x$ is generated randomly with
a uniform distribution in $(0,1)$, and each experiment is repeated $10$ times.
\begin{figure}[htbp]\label{fig1}
    \centering
  \subfloat[The optimality residual versus iteration
counter]{%
       \includegraphics[width=0.9\linewidth]{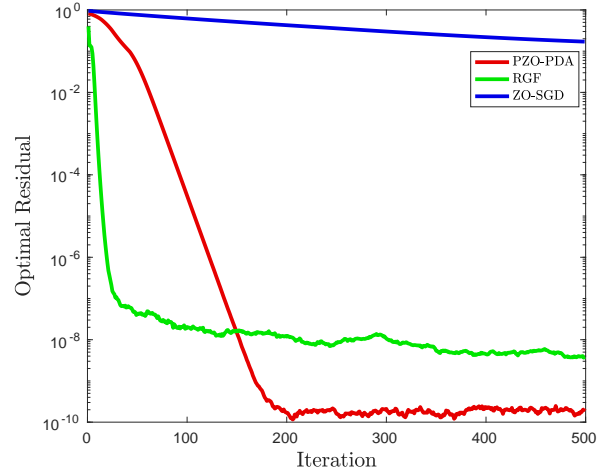}}
    \label{1a}\\
  \subfloat[The constraint violation versus iteration counter]{%
        \includegraphics[width=0.9\linewidth]{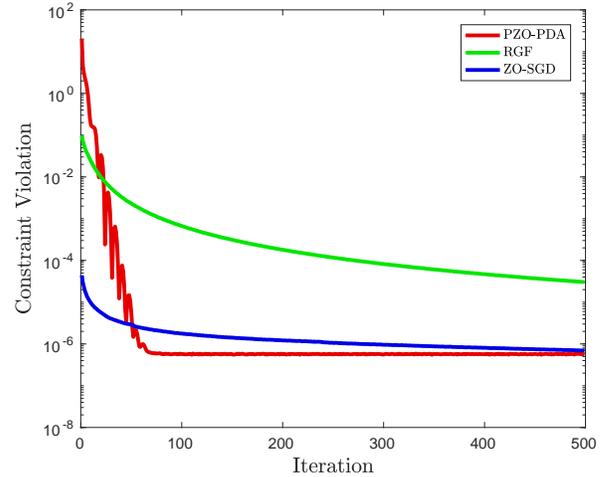}}
    \label{1c}
  \caption{Comparison of different zeroth-order algorithms on nonconvex PCA problem with {$\alpha=10^{-4}$}.}
  \label{fig1} 
\end{figure}
The numerical results are illustrated in Fig. \ref{fig1}. We compare PZO-PDA with RGF algorithm \cite{nesterov2011random} with diminishing step size $0.01\sqrt{\log(2)}/r$, and ZO-SGD algorithm \cite{ghadimi2013stochastic} using decreasing step size $0.01/\sqrt{r}$.
For the performance metric, we use the optimal residual $\norm{x-\Po_{h+\iota_X}\left[x-\nabla f(x)\right]}^2$, and constraint violation.

PZO-PDA exhibits fast convergence to the optimal
solution with even proper constant step size compared to other zeroth-order type algorithms. Although the optimal residual for RGF vanishes comparably to the new algorithm, but it violates linear constraints significantly in contrast to other zeroth-order algorithms. 
%\section*{ACKNOWLEDGMENT}
%This work was partially supported by NSF grant IIS-1741431.

%\bibliography{main}
\bibliography{main}
%\printbibliography
\appendix
\begin{appendices}
\subsection{Proof of Lemma \ref{Titerationdecay}}
\proof
Since ${\beta > 3L_\mu+1}$,  by changing of $x$ we obtain for $C$,
\begin{gather*}
\begin{split}
& C(x^{r+1},\lambda^{r}) - C(x^{r},\lambda^{r}) \\
&= C(x^{r+1},\lambda^{r})+\frac{\beta}{2}\norm{x^{r+1}-x^r}_{B^TB}^2 - C(x^{r},\lambda^{r})\\
&~~~ -\frac{\beta}{2}\norm{x^{r+1}-x^r}_{B^TB}^2\\
& \stackrel{a}{\leq} \left<\nabla f_{\mu}(x^{r+1})+\eta^{r+1} + (1-\rho\gamma)A^T\lambda^r+\rho A^T(Ax^{r+1}-b)\right.\\
&\qquad  \left. +\beta  B^TB(x^{r+1}-x^{r}), x^{r+1}-x^{r}\right>-\frac{\beta-L_{\mu}}{2}\norm{x^{r+1}-x^{r}}^2\\
& = \left<\nabla f_{\mu}(x^{r+1})-\overline{G}_{\mu}^{r}+\overline{G}_{\mu}^{r}+\eta^{r+1} + (1-\rho\gamma)A^T\lambda^r\right.\\
&\qquad  \left. +\rho A^T(Ax^{r+1}-b)+\beta B^TB(x^{r+1}-x^{r}),x^{r+1}-x^{r}\right>\\
&~~~-\frac{\beta-L_{\mu}}{2}\norm{x^{r+1}-x^{r}}^2\\
& \stackrel{b}{\leq} \left<\nabla f_{\mu}(x^{r+1})-\overline{G}_{\mu}^{r}, x^{r+1}-x^{r}\right>-\frac{\beta-L_{\mu}}{2}\norm{x^{r+1}-x^{r}}^2\\
& \leq \frac{1}{2}\frac{\tilde{\sigma}^2}{J} + \frac{1}{2}\norm{x^{r+1}-x^{r}}^2-\left(\frac{\beta-3L_{\mu}}{2}\right)\norm{x^{r+1}-x^{r}}^2\\
&  = \frac{1}{2}\frac{\tilde{\sigma}^2}{J} -\left(\frac{\beta-3L_{\mu}-1}{2}\right)\norm{x^{r+1}-x^{r}}^2
\end{split}
\end{gather*}
 where in $\stackrel{a}{\leq}$ we have used the fact that ${\beta > L_\mu}$, $\rho \geq \beta$, $A^TA+B^TB\geq I_M$, and strong convexity of function $C(x,\lambda^{r})+\frac{\beta}{2}\norm{x-x^r}^{2}_{B^TB}$ with modulus $\beta-L_{\mu}$ [here ${\eta^{r+1}}\in\partial h(x^{r+1})$]. Note that $\stackrel{b}{\leq}$ is true due to the optimality condition \eqref{lem1:eq18} for $x$-subproblem. The last inequality is due to the fact that $f_\mu(x)$ is $L_\mu$-smooth. 
By Lemma 2 of \cite{hajinezhadperturbed}, we further have
  \begin{equation}
 \begin{split}
 &C(x^{r+1},\lambda^{r+1})-C(x^{r+1},\lambda^{r})\\
  &~~~~= (1-\rho\gamma) \left(\frac{1}{\rho}\norm{\lambda^{r+1}-\lambda^{r}}^2\right.\\
  & \left.~~~~~~+ \frac{\gamma}{2}(\norm{\lambda^{r}}^2 -\norm{\lambda^{r+1}}^2 -\norm{\lambda^{r+1}-\lambda^{r}}^2)\right). 
 \end{split}
 \end{equation}
 Summarizing the above arguments, we obtain the desired inequality in \eqref{lem2:25}.
\endproof
\subsection{Proof of Lemma \ref{lemma:potential:lowbd}}
\proof
We have the following 
\begin{gather*}
\begin{align*}
&Q(x^{r+1},\lambda^{r+1};x^{r},\lambda^{r}) - Q(x^{r},\lambda^{r};x^{r-1},\lambda^{r-1})\\
& = C(x^{r+1},\lambda^{r+1}) + \frac{(1-\rho\gamma)\gamma}{2}\norm{\lambda^{r+1}}^2\\
& ~~~~+\frac{1}{2}\left(\frac{(1-\rho\gamma)}{\rho}\norm{\lambda^{r+1}-\lambda^{r}}^2\right.\\
&\left.~~~~~~~~~~+\beta \norm{x^{r+1}-x^{r}}^2_{B^TB}+3L_{\mu}^2\norm{x^{r+1}-x^{r}}^2\right)\\
&~~~~- C(x^{r},\lambda^{r}) - \frac{(1-\rho\gamma)\gamma}{2}\norm{\lambda^{r}}^2\\
&~~~~-\frac{1}{2}\left(\frac{(1-\rho\gamma)}{\rho}\norm{\lambda^{r}-\lambda^{r-1}}^2\right.\\
&~~~~\left.+\beta \norm{x^{r}-x^{r-1}}^2_{B^TB}+3L_{\mu}^2\norm{x^{r}-x^{r-1}}^2\right)\\
&\stackrel{a}{\leq} - \frac{(1-\gamma\rho)\gamma}{2}\norm{\lambda^{r+1}}^2 + \frac{(1-\gamma\rho)\gamma}{2}\norm{\lambda^r}^2\\
&~~~~+\left(\frac{(1-\rho\gamma)(2-\rho\gamma)}{2\rho}\right)\norm{\lambda^{r+1}-\lambda^{r}}^2\\
& ~~~~+ \frac{1}{2}\frac{\tilde{\sigma}^2}{J} -\left(\frac{\beta-3L_{\mu}-1}{2}\right)\norm{x^{r+1}-x^{r}}^2\\
& ~~~~+\frac{(1-\gamma\rho)\gamma}{2}\norm{\lambda^{r+1}}^2 - \frac{(1-\gamma\rho)\gamma}{2}\norm{\lambda^r}^2\\
& ~~~~+\frac{1}{2}\left(\frac{(1-\rho\gamma)}{\rho}\norm{\lambda^{r+1}-\lambda^{r}}^2\right.\\
&\left.~~~~+\beta \norm{x^{r+1}-x^{r}}^2_{B^TB}+3L_{\mu}^2\norm{x^{r+1}-x^{r}}^2\right)\\
& ~~~~-\frac{1}{2}\left(\frac{(1-\rho\gamma)}{\rho}\norm{\lambda^{r}-\lambda^{r-1}}^2\right.\\
&\left.~~~~+\beta \norm{x^{r}-x^{r-1}}^2_{B^TB}+3L_{\mu}^2\norm{x^{r}-x^{r-1}}^2\right)\\
\end{align*}
\end{gather*}
where in $\stackrel{a}{\leq}$ we used Lemma \ref{Titerationdecay}. Therefore, from the above inequality we obtain
\begin{align*}
& Q(x^{r+1},\lambda^{r+1};x^{r},\lambda^{r}) - Q(x^{r},\lambda^{r};x^{r-1},\lambda^{r-1})\\
&\stackrel{a}{\leq}\left(\frac{(1-\rho\gamma)(2-\rho\gamma)}{2\rho}\right)\norm{\lambda^{r+1}-\lambda^{r}}^2\\
&~~~~+ \frac{1}{2}\frac{\tilde{\sigma}^2}{J} -\left(\frac{\beta-3L_{\mu}-1}{2}\right)\norm{x^{r+1}-x^{r}}^2\\
&~~~~+\frac{(1-\rho\gamma)}{2\rho} \norm{\lambda^{r}-\lambda^{r-1}}^2+\frac{\beta}{2} \norm{x^{r}-x^{r-1}}_{B^TB}^2\\
&~~~~+{3}\frac{\tilde{\sigma}^2}{J}+\frac{3L_{\mu}^2}{2}\norm{x^{r}-x^{r-1}}^2\\
&~~~~+ \frac{1}{2}\norm{x^{r+1}-x^{r}}^2 -\gamma\norm{\lambda^{r+1}-\lambda^{r}}^2 + \frac{3L_{\mu}^2}{2}\norm{x^{r+1}-x^{r}}^2\\
&~~~~-\frac{1}{2}\left(\frac{(1-\rho\gamma)}{\rho}\norm{\lambda^{r}-\lambda^{r-1}}^2\right.\\
&~~~~\left.+\beta \norm{x^{r}-x^{r-1}}^2_{B^TB}+3L_{\mu}^2\norm{x^{r}-x^{r-1}}^2\right)\\
&= -\left(\frac{(1-\rho\gamma)(\rho\gamma-2)}{2\rho}+\gamma\right)\norm{\lambda^{r+1}-\lambda^{r}}^2\\
&~~~~-\left(\frac{\beta-1}{2}-\frac{3 L_{\mu}^2}{2}-\frac{3L_{\mu}}{2}-\frac{1}{2}\right)\norm{x^{r+1}-x^{r}}^2 + \frac{7}{2}\frac{\tilde{\sigma}^2}{J} \\
& = -a_1\norm{\lambda^{r+1}-\lambda^{r}}^2 -a_2\norm{x^{r+1}-x^{r}}^2 + \frac{7}{2}\frac{\tilde{\sigma}^2}{J}
\end{align*}
with $ a_1=\left(\frac{(1-\rho\gamma){\gamma}}{2}+\gamma-\frac{1-\rho\gamma}{\rho}\right)$, $ a_2=\left(\frac{\beta}{2}-\frac{3 L_{\mu}^2}{2}-\frac{3L_{\mu}}{2}-1\right)$,  where $\stackrel{a}{\leq}$ implied by Lemma \ref{iterationdecay}.
Therefore, in order to make the potential function decrease, it is sufficient to have 
\begin{equation*}
\frac{(1-\rho\gamma){\gamma}}{2} + \gamma - \frac{1-\rho\gamma}{\rho} >0,\,\,\text{and}\, \beta > (3+3 L_\mu)L_\mu+2.
\end{equation*} 
\endproof
\end{appendices}
\end{document}